\documentstyle[leqno]{article}

\newtheorem{th}{Theorem}[section]
\newtheorem{prop}[th]{Proposition}
\newtheorem{cor}[th]{Corollary}
\newcounter{defin}[section]
\renewcommand{\thedefin}{\thesection.\arabic{defin}}
\newcounter{ex}[section]
\renewcommand{\theex}{\thesection.\arabic{ex}}
\newcounter{rem}[section]
\renewcommand{\therem}{\thesection.\arabic{rem}}
\title{Torsions and Curvatures \\
on Jet Fibre Bundle $J^1(T,M)$}
\date{}
\author{Mircea Neagu and Constantin Udri\c ste}
\begin{document}
\maketitle
\begin{abstract}
The aim of this paper is twofold. On the one hand, to study the local representations
of d-connections, d-torsions, and d-curvatures with respect to an adapted basis
on the jet fibre bundle of order one. On the other  hand, to open the
problem of prolongations of tensors and connections from  a product of two
manifolds to 1-jet fibre bundle associated to these manifolds.
Section 1 defines the notion of $\Gamma$-linear connection on the jet fibre
bundle of order one and determines its nine local  components.
Section 2 studies the main twelve components of torsion d-tensor field, and
Section 3 describes the eighteen components of curvature d-tensor  field. Via
the Ricci and Bianchi identities,
Section 4 emphasizes the non-independence of the torsion and curvature d-tensors.
Section 5 studies the problem of prolongation of vector fields from $T\times M$ to
1-jet space $J^1(T,M)$.
\end{abstract}

\section{Components of $\Gamma$-linear connections}

\hspace{5mm}  Let us consider $T$ (resp. $M$) a "{\it temporal} " (resp.
"{\it spatial} ") manifold of dimension $p$ (resp. $n$), coordinated by
$(t^\alpha)_{\alpha=\overline{1,p}}$ (resp. $(x^i)_{i=\overline{1,n}})$. Let
\linebreak $J^1(T,M)\to T\times M$ be the jet fibre bundle of order one
associated to these manifolds. The {\it bundle of configuration} $J^1(T,M)$
is coordinated by $(t^\alpha,x^i,x^i_\alpha)$, where $\alpha=\overline{1,p}$
and $i=\overline{1,n}$ \cite{10}. Note that, throughout this paper, the indices
$\alpha,\beta,\gamma,\ldots$ run from \linebreak 1 to $p$ and the indices $i,j,k\ldots$
run from 1 to $n$.

On $E=J^1(T,M)$, we fixe a nonlinear connection $\Gamma$ defined by the
{\it temporal} components $M^{(i)}_{(\alpha)\beta}$ and the {\it spatial}
components $N^{(i)}_{(\alpha)j}$. We recall that the transformation rules of
the local components of the nonlinear connection $\Gamma$ are expressed
by \cite{10}
\begin{equation}
\left\{\begin{array}{l}\medskip
\displaystyle{\tilde M^{(j)}_{(\beta)\mu}{\partial\tilde t^\mu\over\partial
t^\alpha}=M^{(k)}_{(\gamma)\alpha}{\partial\tilde x^j\over\partial x^k}
{\partial t^\gamma\over\partial\tilde t^\beta}-{\partial\tilde x^j_\beta\over
\partial t^\alpha}}\\
\displaystyle{\tilde N^{(j)}_{(\beta)k}{\partial\tilde x^k\over\partial
x^i}=N^{(k)}_{(\gamma)i}{\partial\tilde x^j\over\partial x^k}
{\partial t^\gamma\over\partial\tilde t^\beta}-{\partial\tilde x^j_\beta\over
\partial x^i}}.
\end{array}\right.
\end{equation}
Let $\displaystyle{\left\{{\delta\over\delta t^\alpha},{\delta\over\delta x^i},
{\partial\over\partial x^i_\alpha}\right\}\subset{\cal X}(E)}$ and $\{dt^\alpha,
dx^i,\delta x^i_\alpha\}\subset{\cal X}^*(E)$ be the dual adapted bases associated
to the nonlinear connection $\Gamma=(M^{(i)}_{(\alpha)\beta}, N^{(i)}_{(\alpha)
j})$, by the formulas
\begin{equation}
\left\{\begin{array}{l}\medskip
\displaystyle{
{\delta\over\delta t^\alpha}={\partial\over\partial t^\alpha}-
M^{(j)}_{(\beta)\alpha}{\partial\over\partial x^j_\beta}}\\\medskip
\displaystyle{{\delta\over\delta x^i}={\partial\over\partial x^i}-
N^{(j)}_{(\beta)i}{\partial\over\partial x^j_\beta}}\\
\delta x^i_\alpha=dx^i_\alpha+M^{(i)}_{(\alpha)\beta}dt^\beta+N^{(i)}_{(
\alpha)j}dx^j.
\end{array}\right.
\end{equation}
These  bases will be used in the description of geometrical objects on $E$,
because their transformation laws are very simple \cite{10}:
\begin{equation}\label{tlab}
\left\{\begin{array}{lll}\medskip
\displaystyle{{\delta\over\delta t^\alpha}={\partial\tilde t^\beta\over\partial
t^\alpha}{\delta\over\delta\tilde t^\beta}},&
\displaystyle{{\delta\over\delta x^i}={\partial\tilde x^j\over\partial
x^i}{\delta\over\delta\tilde x^j}},&
\displaystyle{{\partial\over\partial x^i_\alpha}={\partial\tilde x^j\over\partial
x^i}{\partial t^\alpha\over\partial\tilde t^\beta}{\partial\over\partial\tilde
x^j_\beta}},\\
\displaystyle{dt^\alpha={\partial t^\alpha\over\partial\tilde
t^\beta}d\tilde t^\beta},&
\displaystyle{dx^i={\partial x^i\over\partial\tilde
x^j}d\tilde x^j},&
\displaystyle{\delta x^i_\alpha={\partial x^i\over\partial\tilde
x^j}{\partial\tilde t^\beta\over\partial t^\alpha}\delta\tilde
x^j_\beta}.
\end{array}\right.
\end{equation}

In order to develope the theory of $\Gamma$-linear connections on the 1-jet
space $E$, we need the following
\begin{prop}
i) The Lie algebra ${\cal X}(E)$ of vector fields decomposes as
$$
{\cal X}(E)={\cal X}({\cal H}_T)\oplus{\cal X}({\cal H}_M)\oplus
{\cal X}({\cal V}),
$$
where
$$
{\cal X}({\cal H}_T)=Span\left\{{\delta\over\delta t^\alpha}\right\},\quad
{\cal X}({\cal H}_M)=Span\left\{{\delta\over\delta x^i}\right\},\quad
{\cal X}({\cal V})=Span\left\{{\partial\over\partial x^i_\alpha}\right\}.
$$

ii) The Lie algebra ${\cal X}^*(E)$ of covector fields decomposes as
$$
{\cal X}^*(E)={\cal X}^*({\cal H}_T)\oplus{\cal X}^*({\cal H}_M)\oplus
{\cal X}^*({\cal V}),
$$
where
$$
{\cal X}^*({\cal H}_T)=Span\{dt^\alpha\},\quad
{\cal X}^*({\cal H}_M)=Span\{dx^i\},\quad
{\cal X}^*({\cal V})=Span\{\delta x^i_\alpha\}.
$$
\end{prop}

Let us consider $h_T$, $h_M$ (horizontal) and $v$ (vertical) as the canonical projections of the above
decompositions. In this context, we have
\begin{cor}
i) Any vector field $X$ can be written in the form
$$
X=h_TX+h_MX+vX.
$$

ii) Any covector field $\omega$ can be written in the form
$$
\omega=h_T\omega+h_M\omega+v\omega .
$$
\end{cor}
\addtocounter{defin}{1}
{\bf Definition \thedefin} A linear connection $\nabla:{\cal X}(E)\times{\cal X}(E)\to
{\cal X}(E)$ is called {\it a $\Gamma$-linear connection on $E$} if $\nabla
h_T=0$, $\nabla h_M=0$ and $\nabla v=0$.

In order to describe in local terms a $\Gamma$-linear connection $\nabla$ on
$E$, we need nine unique local components,
\begin{equation}\label{lglc}
\nabla\Gamma=(\bar G^\alpha_{\beta\gamma},\;G^k_{i\gamma},\;G^{(i)(\beta)}_
{(\alpha)(j)\gamma},\;\bar L^\alpha_{\beta j},\;L^k_{ij},\;L^{(i)(\beta)}_{(\alpha)
(j)k},\;\bar C^{\alpha(\gamma)}_{\beta(k)},\;C^{j(\gamma)}_{i(k)},\;C^{(i)(\beta)
(\gamma)}_{(\alpha)(j)(k)}),
\end{equation}
which are locally defined by the relations
$$\begin{array}{l}\medskip
(h_T)\hspace{4mm}\displaystyle{
\nabla_{{\delta\over\delta t^\gamma}}{\delta\over\delta t^\beta}=
\bar G^\alpha_{\beta\gamma}{\delta\over\delta t^\alpha},\;\;
\nabla_{{\delta\over\delta t^\gamma}}{\delta\over\delta x^i}=G^k_{i\gamma}
{\delta\over\delta x^k},\;\;\nabla_{{\delta\over\delta t^\gamma}}
{\partial\over\partial x^i_\beta}=G^{(k)(\beta)}_{(\alpha)(i)\gamma}
{\partial\over\partial x^k_\alpha}},\\\medskip
(h_M)\hspace{3mm}\displaystyle{\nabla_{{\delta\over\delta x^j}}
{\delta\over\delta t^\beta}=\bar L^\alpha_{\beta j}{\delta\over\delta t^
\alpha},\;\;\nabla_{{\delta\over\delta x^j}}{\delta\over\delta x^i}=L^k_{ij}
{\delta\over\delta x^k},\;\;\nabla_{{\delta\over\delta x^j}}
{\partial\over\partial x^i_\beta}=L^{(k)(\beta)}_{(\alpha)(i)j}
{\partial\over\partial x^k_\alpha}},\\
(v)\hspace{6mm}\displaystyle{\nabla_{{\partial\over\partial x^j_\gamma}}
{\delta\over\delta t^\beta}=\bar C^{\alpha(\gamma)}_{\beta(j)}
{\delta\over\delta t^\alpha},\;\nabla_{{\partial\over\partial x^j_\gamma}}
{\delta\over\delta x^i}=C^{k(\gamma)}_{i(j)}
{\delta\over\delta x^k},\;\nabla_{{\partial\over\partial x^j_\gamma}}
{\partial\over\partial x^i_\beta}=C^{(k)(\beta)(\gamma)}_{(\alpha)(i)(j)}
{\partial\over\partial x^k_\alpha}}.
\end{array}
$$

The transformation laws of the elements
$\displaystyle{\left\{{\delta\over\delta t^\alpha},{\delta\over\delta x^i},
{\partial\over\partial x^i_\alpha}\right\}}$ together with
the properties of the $\Gamma$-linear connection $\nabla$, imply
\begin{th}\label{glc}
i) The components of the $\Gamma$-linear connection $\nabla$ modify
by the rules\\\\
$(h_T)\hspace{6mm}
\left\{\begin{array}{l}\medskip\displaystyle{
\bar G^\delta_{\alpha\beta}{\partial\tilde t^\varepsilon\over\partial t^\delta}
=\tilde{\bar G}^\varepsilon_{\mu\gamma}{\partial\tilde t^\mu\over\partial
t^\alpha}{\partial\tilde t^\gamma\over\partial t^\beta}+{\partial^2\tilde t^
\varepsilon\over\partial t^\alpha\partial t^\beta}}\\\medskip
\displaystyle{G^k_{i\gamma}=\tilde G^m_{j\beta}{\partial x^k\over\partial
\tilde x^m}{\partial\tilde x^j\over\partial x^i}{\partial\tilde t^\beta\over
\partial t^\gamma}}\\\medskip
\displaystyle{
G^{(k)(\beta)}_{(\gamma)(i)\alpha}=\tilde G^{(p)(\eta)}_
{(\varepsilon)(j)\mu}{\partial x^k\over\partial\tilde x^p}
{\partial\tilde t^\varepsilon\over\partial t^\gamma}
{\partial\tilde x^j\over\partial x^i}
{\partial t^\beta\over\partial\tilde t^\eta}
{\partial\tilde t^\mu\over\partial t^\alpha}+\delta^k_i
{\partial\tilde t^\mu\over\partial t^\alpha}
{\partial\tilde t^\varepsilon\over\partial t^\gamma}
{\partial^2 t^\beta\over\partial\tilde t^\mu\partial\tilde t^\varepsilon}},
\end{array}\right.\medskip
$\\
$(h_M)\hspace{5mm}
\left\{\begin{array}{l}\medskip\displaystyle{
\bar L^\gamma_{\beta j}{\partial x^j\over\partial\tilde x^l}
=\tilde{\bar L}^\eta_{\mu l}{\partial t^\gamma\over\partial\tilde
t^\eta}{\partial\tilde t^\mu\over\partial t^\beta}}\\\medskip
\displaystyle{L^m_{ij}{\partial\tilde x^r\over\partial x^m}=\tilde L^r_{pq}
{\partial\tilde x^p\over\partial x^i}{\partial\tilde x^q\over\partial x^j}+
{\partial^2\tilde x^r\over\partial x^i\partial x^j}}\\\medskip
\displaystyle{L^{(k)(\beta)}_{(\gamma)(i)j}=\tilde L^{(r)(\eta)}_
{(\nu)(p)l}{\partial x^k\over\partial\tilde x^r}
{\partial\tilde t^\nu\over\partial t^\gamma}
{\partial\tilde x^p\over\partial x^i}
{\partial t^\beta\over\partial\tilde t^\eta}
{\partial\tilde x^l\over\partial x^j}+\delta^\beta_\gamma
{\partial x^k\over\partial\tilde x^r}
{\partial^2\tilde x^r\over\partial x^i\partial x^j}},
\end{array}\right.\medskip
$\\
$(v)\hspace{8mm}
\left\{\begin{array}{l}\medskip\displaystyle{
\bar C^{\gamma(\alpha)}_{\beta(i)}=\tilde{\bar C}^{\mu(\delta)}_
{\varepsilon(j)}{\partial t^\gamma\over\partial\tilde t^\mu}
{\partial\tilde t^\varepsilon\over\partial t^\beta}
{\partial\tilde x^j\over\partial x^i}
{\partial t^\alpha\over\partial\tilde t^\delta}}\\\medskip
\displaystyle{C^{k(\alpha)}_{i(j)}=\tilde C^{s(\beta)}_
{p(r)}{\partial x^k\over\partial\tilde x^s}
{\partial\tilde x^p\over\partial x^i}
{\partial\tilde x^r\over\partial x^j}
{\partial t^\alpha\over\partial\tilde t^\beta}}\\\medskip
\displaystyle{C^{(k)(\beta)(\alpha)}_{(\gamma)(i)(j)}=\tilde C^{(r)(\mu)(\nu)}
_{(\varepsilon)(p)(q)}{\partial x^k\over\partial\tilde x^r}
{\partial\tilde t^\varepsilon\over\partial t^\gamma}
{\partial\tilde x^p\over\partial x^i}
{\partial t^\beta\over\partial\tilde t^\mu}
{\partial\tilde x^q\over\partial x^j}
{\partial t^\alpha\over\partial\tilde t^\nu}}.
\end{array}\right.\medskip
$

ii) Conversely, to give a $\Gamma$-linear connection $\nabla$ on the 1-jet space
$E$ is equivalent to give a set of nine local components \ref{lglc},
whose local transformations laws are described in i.
\end{th}

Theorem \ref{glc} allows us to construct on the 1-jet space $E$ a natural
example of $\Gamma$-linear connection.\medskip\\
\addtocounter{ex}{1}
{\bf Example \theex} Suppose that $h_{\alpha\beta}(t)$ (resp. $\varphi_{ij}(x)$) is
a pseudo-Riemannian metric on $T$ (resp. $M$). We denote $H^\gamma_{\alpha
\beta}$ (resp. $\gamma^k_{ij}$) the Christoffel symbols of the metric $h_{\alpha
\beta}$ (resp. $\varphi_{ij}$). The {\it canonical nonlinear
connection $\Gamma_0$ associated to these metrics} is defined by the local
components \cite{10}
\begin{equation}\label{cnc}
M^{(i)}_{(\alpha)\beta}=-H^\gamma_{\alpha\beta}x^i_\gamma,\quad
N^{(i)}_{(\alpha)j}=\gamma^i_{jm}x^m_\alpha.
\end{equation}
In this context, using the well known local transformation rules of the
Christoffel symbols $H^\gamma_{\alpha\beta}$ and $\gamma^i_{jk}$ and setting
\begin{equation}\label{bglc}
\bar G^\gamma_{\alpha\beta}=H^\gamma_{\alpha\beta},\;\;
G^{(k)(\beta)}_{(\gamma)(i)\alpha}=-\delta^k_iH^\beta_{\alpha\gamma},\;\;
L^k_{ij}=\gamma^k_{ij},\;\;
L^{(k)(\beta)}_{(\gamma)(i)j}=\delta^\beta_\gamma\gamma^k_{ij},
\end{equation}
we conclude that the set of local components
$$
B\Gamma_0=(H^\gamma_{\alpha\beta},\;0,\;G^{(k)(\beta)}_{(\gamma)(i)\alpha},
\;0,\;\gamma^k_{ij},\;L^{(k)(\beta)}_{(\gamma)(i)j},\;0,\;0,\;0)
$$
defines a $\Gamma_0$-linear connection. This is called the {\it Berwald connection
attached to the metrics pair $(h_{\alpha\beta},\varphi_{ij})$.}\medskip\\
\addtocounter{rem}{1}
{\bf Remark \therem} In the particular case $(T,h)=(R,\delta)$, the Berwald
connection reduces to that naturally induced by the canonical spray
$2G^i=\gamma^i_{jk}y^jy^k$ of the classical theory of Lagrange spaces . For more
details, see \cite{5}, \cite{10}.\medskip

Now, let $\nabla$ be a $\Gamma$-linear connection on $E$, locally defined by
\ref{lglc}.
The linear connection $\nabla$ induces a natural linear connection on the
d-tensors set of the jet fibre bundle $E=J^1(T,M)$, in the following fashion:
starting with a vector field $X$ and a d-tensor field $D$
locally expressed by
$$\begin{array}{l}\medskip\displaystyle{
X=X^\alpha{\delta\over\delta t^\alpha}+X^m{\delta\over\delta x^m}
+X^{(m)}_{(\alpha)}{\partial\over\partial x^m_\alpha},}\\
\displaystyle{
D=D^{\alpha i(j)(\delta)\ldots}_{\gamma k(\beta)(l)\ldots}{\delta\over\delta
t^\alpha}\otimes{\delta\over\delta x^i}\otimes{\partial\over\partial x^j_
\beta}\otimes dt^\gamma\otimes dx^k\otimes\delta x^l_\delta\ldots,}
\end{array}
$$
we introduce the covariant derivative
$$
\begin{array}{l}\medskip
\nabla_XD=X^\varepsilon\nabla_{\delta\over\delta t^\varepsilon}D+X^p\nabla_
{\delta\over\delta x^p}D+X^{(p)}_{(\varepsilon)}\nabla_{\partial\over\partial
x^p_\varepsilon}D=\left\{X^\varepsilon D^{\alpha i(j)(\delta)\ldots}_{\gamma k
(\beta)(l)\ldots /\varepsilon}+X^p\right.\\\left.
D^{\alpha i(j)(\delta)\ldots}_{\gamma k(\beta)(l)\ldots\vert p}
+X^{(p)}_{(\varepsilon)}D^{\alpha i(j)(\delta)\ldots}_{\gamma k (\beta)(l)
\ldots}\vert_{(p)}^{(\varepsilon)}\right\}\displaystyle{
{\delta\over\delta t^\alpha}\otimes{\delta\over\delta x^i}\otimes{\partial
\over\partial x^j_\beta}\otimes dt^\gamma\otimes dx^k\otimes\delta x^l_\delta
\ldots,}
\end{array}
$$
where\\\\
$
(h_T)\hspace{6mm}\left\{\begin{array}{l}\medskip\displaystyle{
D^{\alpha i(j)(\delta)\ldots}_{\gamma k(\beta)(l)\ldots /\varepsilon}=
{\delta D^{\alpha i(j)(\delta)\ldots}_{\gamma k(\beta)(l)\ldots}\over
\delta t^\varepsilon}+D^{\mu i(j)(\delta)\ldots}_{\gamma k(\beta)(l)\ldots}
\bar G^\alpha_{\mu\varepsilon}+}\\\medskip
D^{\alpha m(j)(\delta)\ldots}_{\gamma k(\beta)(l)\ldots}G^i_{m\varepsilon}+
D^{\alpha i(m)(\delta)\ldots}_{\gamma k(\mu)(l)\ldots}G^{(j)(\mu)}_{(\beta)
(m)\varepsilon}+\ldots-\\
-D^{\alpha i(j)(\delta)\ldots}_{\mu k(\beta)(l)\ldots}\bar G^\mu_{\gamma
\varepsilon}-D^{\alpha i(j)(\delta)\ldots}_{\gamma m(\beta)(l)\ldots} G^m_
{k\varepsilon}-D^{\alpha i(j)(\mu)\ldots}_{\gamma k(\beta)(m)\ldots} G^{(m)
(\delta)}_{(\mu)(l)\varepsilon}-\ldots,
\end{array}\right.\medskip
$
$
(h_M)\hspace{5mm}\left\{\begin{array}{l}\medskip\displaystyle{
D^{\alpha i(j)(\delta)\ldots}_{\gamma k(\beta)(l)\ldots\vert p}=
{\delta D^{\alpha i(j)(\delta)\ldots}_{\gamma k(\beta)(l)\ldots}\over
\delta x^p}+D^{\mu i(j)(\delta)\ldots}_{\gamma k(\beta)(l)\ldots}
\bar L^\alpha_{\mu p}+}\\\medskip
D^{\alpha m(j)(\delta)\ldots}_{\gamma k(\beta)(l)\ldots}L^i_{mp}+
D^{\alpha i(m)(\delta)\ldots}_{\gamma k(\mu)(l)\ldots}L^{(j)(\mu)}_{(\beta)
(m)p}+\ldots-\\
-D^{\alpha i(j)(\delta)\ldots}_{\mu k(\beta)(l)\ldots}\bar L^\mu_{\gamma
p}-D^{\alpha i(j)(\delta)\ldots}_{\gamma m(\beta)(l)\ldots} L^m_
{kp}-D^{\alpha i(j)(\mu)\ldots}_{\gamma k(\beta)(m)\ldots} L^{(m)
(\delta)}_{(\mu)(l)p}-\ldots,
\end{array}\right.\medskip
$
$
(v)\hspace{8mm}\left\{\begin{array}{l}\medskip\displaystyle{
D^{\alpha i(j)(\delta)\ldots}_{\gamma k(\beta)(l)\ldots}\vert_{(p)}^{(\varepsilon)}=
{\partial D^{\alpha i(j)(\delta)\ldots}_{\gamma k(\beta)(l)\ldots}\over
\partial x^p_\varepsilon}+D^{\mu i(j)(\delta)\ldots}_{\gamma k(\beta)(l)\ldots}
\bar C^{\alpha(\varepsilon)}_{\mu(p)}+}\\\medskip
D^{\alpha m(j)(\delta)\ldots}_{\gamma k(\beta)(l)\ldots}C^{i(\varepsilon)}_{m(p)}+
D^{\alpha i(m)(\delta)\ldots}_{\gamma k(\mu)(l)\ldots}C^{(j)(\mu)(\varepsilon)}
_{(\beta)(m)(p)}+\ldots-\\
-D^{\alpha i(j)(\delta)\ldots}_{\mu k(\beta)(l)\ldots}\bar C^{\mu(\varepsilon)}_
{\gamma(p)}-D^{\alpha i(j)(\delta)\ldots}_{\gamma m(\beta)(l)\ldots} C^{m(\varepsilon)}
_{k(p)}-D^{\alpha i(j)(\mu)\ldots}_{\gamma k(\beta)(m)\ldots} C^{(m)(\delta)
(\varepsilon)}_{(\mu)(l)(p)}-\ldots.
\end{array}\right.\medskip
$

The local operators "$ _{/\varepsilon}$", "$_{\vert p}$" and "$\vert
^{(\varepsilon)}_{(p)}$" are called the {\it $T$-horizontal covariant derivative,
$M$-horizontal covariant derivative} and {vertical covariant derivative}
of the $\Gamma$ linear connection $\nabla$.\medskip\\
\addtocounter{rem}{1}
{\bf Remarks \therem} i) In the particular case of a function $f(t^\gamma,x^k,x^k_
\gamma)$ on $J^1(T,M)$, the above covariant derivatives reduce to
\begin{equation}
\left\{\begin{array}{l}\medskip
f_{/\varepsilon}=\displaystyle{{\delta f\over\delta t^\varepsilon}={\partial f
\over\partial t^\varepsilon}-M^{(k)}_{(\gamma)\varepsilon}{\partial f\over
\partial x^k_\gamma}}\\\medskip
f_{\vert p}=\displaystyle{{\delta f\over\delta x^p}={\partial f
\over\partial x^p}-N^{(k)}_{(\gamma)p}{\partial f\over\partial x^k_\gamma}}\\
f\vert^{(\varepsilon)}_{(p)}=\displaystyle{{\partial f\over\partial x^p_
\varepsilon}}.
\end{array}\right.
\end{equation}

ii) Particularly, starting with a d-vector field $X$ on $J^1(T,M)$, locally
expressed by
$$
X=X^\alpha{\delta\over\delta t^\alpha}+X^i{\delta\over\delta x^i}+
X^{(i)}_{(\alpha)}{\partial\over\partial x^i_\alpha},
$$
the following expressions of above covariant derivatives hold good:\medskip\\
$(h_T)\hspace*{6mm}
\left\{\begin{array}{l}\medskip
\displaystyle{X^\alpha_{/\varepsilon}={\delta X^\alpha\over\delta t^\varepsilon}
+X^\mu\bar G^\alpha_{\mu\varepsilon}}\\\medskip
\displaystyle{X^i_{/\varepsilon}={\delta X^i\over\delta t^\varepsilon}
+X^mG^i_{m\varepsilon}}\\
\displaystyle{X^{(i)}_{(\alpha)/\varepsilon}={\delta X^{(i)}_{(\alpha)}\over
\delta t^\varepsilon}+X^{(m)}_{(\mu)}G^{(i)(\mu)}_{(\alpha)(m)\varepsilon}},
\end{array}\right.\medskip
$
$(h_M)\hspace*{5mm}
\left\{\begin{array}{l}\medskip
\displaystyle{X^\alpha_{\vert p}={\delta X^\alpha\over\delta x^p}
+X^\mu\bar L^\alpha_{\mu p}}\\\medskip
\displaystyle{X^i_{\vert p}={\delta X^i\over\delta x^p}
+X^mL^i_{mp}}\\
\displaystyle{X^{(i)}_{(\alpha)\vert p}={\delta X^{(i)}_{(\alpha)}\over
\delta x^p}+X^{(m)}_{(\mu)}L^{(i)(\mu)}_{(\alpha)(m)p}},
\end{array}\right.\medskip
$
$(v)\hspace*{8mm}
\left\{\begin{array}{l}\medskip
\displaystyle{X^\alpha\vert^{(\varepsilon)}_{(p)}={\partial X^\alpha\over
\partial x^p_\varepsilon}+X^\mu\bar C^{\alpha(\varepsilon)}_{\mu(p)}}\\\medskip
\displaystyle{X^i\vert^{(\varepsilon)}_{(p)}={\partial X^i\over\partial x^p_
\varepsilon}+X^mC^{i(\varepsilon)}_{m(p)}}\\
\displaystyle{X^{(i)}_{(\alpha)}\vert^{(\varepsilon)}_{(p)}={\partial X^{(i)}_
{(\alpha)}\over\partial x^p_\varepsilon}+X^{(m)}_{(\mu)}C^{(i)(\mu)(\varepsilon)}
_{(\alpha)(m)(p)}}.
\end{array}\right.
$

iii) The local covariant derivatives associated to the Berwald $\Gamma_0$-linear
connection, will be denoted by
"$ _{/\!/\varepsilon}$", "$_{\Vert p}$" and "$\Vert^{(\varepsilon)}_{(p)}$".

Denoting by "$_{:A}$" one of the covariant  derivatives "$_{/\varepsilon}$",
"$_{\vert p}$" or "$\vert^{(\varepsilon)}_{(p)}$", one easily proves the
following

\begin{prop}\label{propr}
If $D^{...}_{...}$ and $F^{...}_{...}$
are two d-tensor fields on $E$, then the following statements hold good.

i) $D^{\ldots}_{\ldots :A}$ are the components of a new d-tensor field,

ii) $(D^{\ldots}_{\ldots}+F^{\ldots}_{\ldots})_{:A}=D^{\ldots}_{\ldots :A}+
F^{\ldots}_{\ldots :A}$,

iii) $(D^{\ldots}_{\ldots}\otimes F^{\ldots}_{\ldots})_{:A}=D^{\ldots}_
{\ldots :A}\otimes F^{\ldots}_{\ldots}+D^{\ldots}_{\ldots}\otimes F^{\ldots}_
{\ldots :A}$,

iv) The operator "$_{:A}$" commutes with the operation of contraction.
\end{prop}

\section{Components of torsion d-tensor field}

\setcounter{equation}{0}
\hspace{5mm} Let us start with a fixed $\Gamma$-linear connection $\nabla$ on
$E=J^1(T,M)$, defined by the local components \ref{lglc}. The torsion d-tensor field associated to $\nabla$ is
$$\mbox{\bf T}:{\cal X}(E)\times{\cal X}(E)\to {\cal X}(E),\quad
\mbox{\bf T}(X,Y)=\nabla_XY-\nabla_YX-[X,Y],\;\;\forall\;
X,Y\in{\cal X}(E).$$

To characterize locally the torsion d-tensor {\bf T} of the connection $\nabla$,
we need the next
\begin{prop}
The following bracket identities are true,
$$
\begin{array}{ll}\medskip
\displaystyle{\left[{\delta\over\delta t^\alpha},{\delta\over\delta t^\beta}
\right]=R^{(m)}_{(\mu)\alpha\beta}{\partial\over\partial x^m_\mu}},&
\displaystyle{\left[{\delta\over\delta t^\alpha},{\delta\over\delta x^j}
\right]=R^{(m)}_{(\mu)\alpha j}{\partial\over\partial x^m_\mu}},\\\medskip
\displaystyle{\left[{\delta\over\delta t^\alpha},{\partial\over\partial x^j_
\beta}\right]={\partial M^{(m)}_{(\mu)\alpha}\over\partial x^j_\beta}
{\partial\over\partial x^m_\mu}},&
\displaystyle{\left[{\delta\over\delta x^i},{\delta\over\delta x^j}
\right]=R^{(m)}_{(\mu)ij}{\partial\over\partial x^m_\mu}},\\
\displaystyle{\left[{\delta\over\delta x^i},{\partial\over\partial x^j_
\beta}\right]={\partial N^{(m)}_{(\mu)i}\over\partial x^j_\beta}
{\partial\over\partial x^m_\mu}},&
\displaystyle{\left[{\partial\over\partial x^i_\alpha},{\partial\over\partial
x^j_\beta}\right]=0},
\end{array}
$$\medskip
where $M^{(m)}_{(\mu)\alpha}$ and $N^{(m)}_{(\mu)i}$ are the local components
of the nonlinear connection $\Gamma$ while the components $R^{(m)}_{(\mu)\alpha\beta},\;
R^{(m)}_{(\mu)\alpha j},\;R^{(m)}_{(\mu)ij}$ are d-tensors expressed by
\medskip\linebreak
$
\displaystyle{
R^{(m)}_{(\mu)\alpha\beta}={\delta M^{(m)}_{(\mu)\alpha}\over
\delta t^\beta}-{\delta M^{(m)}_{(\mu)\beta}\over\delta t^\alpha}},\;\;
\displaystyle{
R^{(m)}_{(\mu)\alpha j}={\delta M^{(m)}_{(\mu)\alpha}\over
\delta x^j}-{\delta N^{(m)}_{(\mu)j}\over\delta t^\alpha},\;\;
R^{(m)}_{(\mu)ij}={\delta N^{(m)}_{(\mu)i}\over
\delta x^j}-{\delta N^{(m)}_{(\mu)j}\over\delta x^i}}.
$
\end{prop}

Consequently, the torsion d-tensor  field {\bf T}
of the $\Gamma$-linear connection $\nabla$ can be described locally by
\begin{th}
The torsion d-tensor {\bf T} of the $\Gamma$-linear connection $\nabla$ is
determined by the following local expressions:
$$
\displaystyle{
h_T\mbox{\bf T}\left({\delta\over\delta t^\beta},{\delta\over\delta
t^\alpha}\right)=\bar T^\mu_{\alpha\beta}{\delta\over\delta t^\mu},\quad
h_M\mbox{\bf T}\left({\delta\over\delta t^\beta},{\delta\over\delta
t^\alpha}\right)=0,}
$$
$$\displaystyle{
v\mbox{\bf T}\left({\delta\over\delta t^\beta},{\delta\over\delta
t^\alpha}\right)=R^{(m)}_{(\mu)\alpha\beta}{\partial\over\partial x^m_\mu}},
$$
$$
\displaystyle{h_T\mbox{\bf T}\left({\delta\over\delta x^j},{\delta\over\delta
t^\alpha}\right)=\bar T^\mu_{\alpha j}{\delta\over\delta t^\mu},\quad
h_M\mbox{\bf T}\left({\delta\over\delta x^j},{\delta\over\delta
t^\alpha}\right)=T^m_{\alpha j}{\delta\over\delta x^m},}
$$
$$
\displaystyle{
v\mbox{\bf T}\left({\delta\over\delta x^j},{\delta\over\delta
t^\alpha}\right)=R^{(m)}_{(\mu)\alpha j}{\partial\over\partial x^m_\mu}},
$$
$$
\displaystyle{h_T\mbox{\bf T}\left({\delta\over\delta x^j},{\delta\over\delta
x^i}\right)=0},\quad
\displaystyle{h_M\mbox{\bf T}\left({\delta\over\delta x^j},{\delta\over\delta
x^i}\right)=T^m_{ij}{\delta\over\delta x^m}},
$$
$$
\displaystyle{v\mbox{\bf T}\left({\delta\over\delta x^j},{\delta\over\delta
x^i}\right)=R^{(m)}_{(\mu)ij}{\partial\over\partial x^m_\mu}},
$$
$$
\displaystyle{
h_T\mbox{\bf T}\left({\partial\over\partial x^j_\beta},{\delta
\over\delta t^\alpha}\right)=\bar P^{\mu(\beta)}_{\alpha(j)}{\delta\over\delta
t^\mu}},\quad
\displaystyle{
h_M\mbox{\bf T}\left({\partial\over\partial x^j_\beta},{\delta
\over\delta t^\alpha}\right)=0},
$$
$$
\displaystyle{
v\mbox{\bf T}\left({\partial\over\partial x^j_\beta},{\delta
\over\delta t^\alpha}\right)=P^{(m)\;(\beta)}_{(\mu)\alpha(j)}{\partial\over
\partial x^m_\mu}},
$$
$$
\displaystyle{
h_T\mbox{\bf T}\left({\partial\over\partial x^j_\beta},{\delta
\over\delta x^i}\right)=0},\quad
\displaystyle{
h_M\mbox{\bf T}\left({\partial\over\partial x^j_\beta},{\delta
\over\delta t^\alpha}\right)=P^{m(\beta)}_{i(j)}{\delta\over\delta x^m}},
$$
$$
\displaystyle{
v\mbox{\bf T}\left({\partial\over\partial x^j_\beta},{\delta
\over\delta x^i}\right)=P^{(m)\;(\beta)}_{(\mu)i(j)}{\partial\over
\partial x^m_\mu}},
$$
$$
\displaystyle{
h_T\mbox{\bf T}\left({\partial\over\partial x^j_\beta},{\partial\over\partial
x^i_\alpha}\right)=0},\quad
\displaystyle{
h_M\mbox{\bf T}\left({\partial\over\partial x^j_\beta},{\partial\over\partial
x^i_\alpha}\right)=0},
$$
$$
\displaystyle{
v\mbox{\bf T}\left({\partial\over\partial x^j_\beta},{\partial\over\partial
x^i_\alpha}\right)=S^{(m)(\alpha)(\beta)}_{(\mu)(i)(j)}{\partial\over
\partial x^m_\mu}},
$$
where $R^{(m)}_{(\mu)\alpha\beta}$, $R^{(m)}_{(\mu)\alpha j}$, $R^{(m)}_{(\mu)ij}$
are the d-tensors constructed above, and\medskip\linebreak
$\bar T^\mu_{\alpha\beta}=\bar G^\mu_{\alpha\beta}-\bar G^\mu_{\beta\alpha}$,
$\bar T^\mu_{\alpha j}=\bar L^\mu_{\alpha j}$,
$\bar P^{\mu(\beta)}_{\alpha(j)}=\bar C^{\mu(\beta)}_{\alpha(j)}$,
$T^m_{\alpha j}=-G^m_{j\alpha}$,\medskip\linebreak
$T^m_{ij}=L^m_{ij}-L^m_{ji}$,
$P^{m(\beta)}_{i(j)}=C^{m(\beta)}_{i(j)}$,
$S^{(m)(\alpha)(\beta)}_{(\mu)(i)(j)}=C^{(m)(\alpha)(\beta)}_{(\mu)(i)(j)}-
C^{(m)(\beta)(\alpha)}_{(\mu)(j)(i)}$,\medskip\linebreak
$\displaystyle{
P^{(m)\;(\beta)}_{(\mu)\alpha(j)}={\partial M^{(m)}_{(\mu)\alpha}\over\partial
x^j_\beta}-G^{(m)(\beta)}_{(\mu)(j)\alpha}},\;
\displaystyle{
P^{(m)\;(\beta)}_{(\mu)i(j)}={\partial N^{(m)}_{(\mu)i}\over\partial
x^j_\beta}-L^{(m)(\beta)}_{(\mu)(j)i}}$.
\end{th}

\begin{cor}
The torsion {\bf\em T} of the $\Gamma$-linear connection $\nabla$ is determined by
twelve effective d-tensor fields, arranged in the following table
\begin{equation}
\begin{tabular}{|c|c|c|c|}
\hline
&$h_T$&$h_M$&$v$\\
\hline
$h_Th_T$&$\bar T^\mu_{\alpha\beta}$&0&$R^{(m)}_{(\mu)\alpha\beta}$\\
\hline
$h_Mh_T$&$\bar T^\mu_{\alpha j}$&$T^m_{\alpha j}$&$R^{(m)}_{(\mu)\alpha j}$\\
\hline
$h_Mh_M$&0&$T^m_{ij}$&$R^{(m)}_{(\mu)ij}$\\
\hline
$vh_T$&$\bar P^{\mu(\beta)}_{\alpha(j)}$&0&$P^{(m)\;(\beta)}_{(\mu)\alpha(j)}$\\
\hline
$vh_M$&0&$P^{m(\beta)}_{i(j)}$&$P^{(m)\;(\beta)}_{(\mu)i(j)}$\\
\hline
$vv$&0&0&$S^{(m)(\alpha)(\beta)}_{(\mu)(i)(j)}$\\
\hline
\end{tabular}
\end{equation}
\end{cor}
\addtocounter{rem}{1}
{\bf Remark \therem} In the particular case of the Berwald $\Gamma_0$-linear connection
associated to the metrics $h_{\alpha\beta}$ and $\varphi_{ij}$, all torsion
d-tensors vanish, except
$$
R^{(m)}_{(\mu)\alpha\beta}=-H^\gamma_{\mu\alpha\beta}x^m_\gamma,\quad
R^{(m)}_{(\mu)ij}=r^m_{ijl}x^l_\mu,
$$
where $H^\gamma_{\mu\beta\gamma}$ (resp. $r^m_{ijl}$) are the curvature
tensors of the metric $h_{\alpha\beta}$ (resp. $\varphi_{ij}$).

\section{Components of curvature d-tensor field}

\setcounter{equation}{0}
\hspace{5mm} From the general theory of linear connections, we recall that the
curvature d-tensor field associated to the $\Gamma$-linear connection
$\nabla$ is
$$
\mbox{\bf R}(X,Y)Z=\nabla_X\nabla_YZ-\nabla_Y\nabla_XZ-\nabla_{[X,Y]}Z,\quad
\forall\;X,Y,Z\in{\cal X}(E).
$$

Using an adapted basis and the properties of the $\Gamma$-linear connection
$\nabla$, one easily prove the following
\begin{th}
The curvature d-tensor {\em\bf R} of the $\Gamma$-linear connection $\nabla$
is determined by the following eighteen local expressions,
$$
\displaystyle{
\mbox{\em\bf R}\left({\delta\over\delta t^\gamma},{\delta\over\delta
t^\beta}\right){\delta\over\delta t^\alpha}=\bar R^\delta_{\alpha\beta\gamma}
{\delta\over\delta t^\delta}},\quad
\displaystyle{
\mbox{\em\bf R}\left({\delta\over\delta t^\gamma},{\delta\over\delta
t^\beta}\right){\delta\over\delta x^i}=R^l_{i\beta\gamma}{\delta
\over\delta x^l}},
$$
$$
\displaystyle{
\mbox{\em\bf R}\left({\delta\over\delta t^\gamma},{\delta\over\delta
t^\beta}\right){\partial\over\partial x^i_\alpha}=R^{(l)(\alpha)}_{(
\delta)(i)\beta\gamma}{\partial\over\partial x^l_\delta}},
$$
$$
\displaystyle{
\mbox{\em\bf R}\left({\delta\over\delta x^k},{\delta\over\delta
t^\beta}\right){\delta\over\delta t^\alpha}=\bar R^\delta_{\alpha\beta k}
{\delta\over\delta t^\delta}},\quad
\displaystyle{
\mbox{\em\bf R}\left({\delta\over\delta x^k},{\delta\over\delta
t^\beta}\right){\delta\over\delta x^i}=R^l_{i\beta k}{\delta
\over\delta x^l}},
$$
$$
\displaystyle{
\mbox{\em\bf R}\left({\delta\over\delta x^k},{\delta\over\delta
t^\beta}\right){\partial\over\partial x^i_\alpha}=R^{(l)(\alpha)}_{(
\delta)(i)\beta k}{\partial\over\partial x^l_\delta}},
$$
$$
\displaystyle{
\mbox{\em\bf R}\left({\delta\over\delta x^k},{\delta\over\delta
x^j}\right){\delta\over\delta t^\alpha}=\bar R^\delta_{\alpha jk}
{\delta\over\delta t^\delta}},\quad
\displaystyle{
\mbox{\em\bf R}\left({\delta\over\delta x^k},{\delta\over\delta
x^j}\right){\delta\over\delta x^i}=R^l_{ijk}{\delta\over\delta x^l}},
$$
$$
\displaystyle{
\mbox{\em\bf R}\left({\delta\over\delta x^k},{\delta\over\delta
x^j}\right){\partial\over\partial x^i_\alpha}=R^{(l)(\alpha)}_{(
\delta)(i)jk}{\partial\over\partial x^l_\delta}},
$$
$$
\displaystyle{
\mbox{\em\bf R}\left({\partial\over\partial x^k_\gamma},{\delta\over\delta
t^\beta}\right){\delta\over\delta t^\alpha}=\bar P^{\delta\;\;(\gamma)}_{\alpha
\beta(k)}{\delta\over\delta t^\delta}},\quad
\displaystyle{
\mbox{\em\bf R}\left({\partial\over\partial x^k_\gamma},{\delta\over\delta
t^\beta}\right){\delta\over\delta x^i}=P^{l\;\;(\gamma)}_{i\beta(k)}{\delta
\over\delta x^l}},
$$
$$
\displaystyle{
\mbox{\em\bf R}\left({\partial\over\partial x^k_\gamma},{\delta\over\delta
t^\beta}\right){\partial\over\partial x^i_\alpha}=P^{(l)(\alpha)\;(\gamma)}_
{(\delta)(i)\beta(k)}{\partial\over\partial x^l_\delta}},
$$
$$
\displaystyle{
\mbox{\em\bf R}\left({\partial\over\partial x^k_\gamma},{\delta\over\delta
x^j}\right){\delta\over\delta t^\alpha}=\bar P^{\delta\;\;(\gamma)}_{\alpha
j(k)}{\delta\over\delta t^\delta}},\quad
\displaystyle{
\mbox{\em\bf R}\left({\partial\over\partial x^k_\gamma},{\delta\over\delta
x^j}\right){\delta\over\delta x^i}=P^{l\;\;(\gamma)}_{ij(k)}{\delta
\over\delta x^l}},
$$
$$
\displaystyle{
\mbox{\em\bf R}\left({\partial\over\partial x^k_\gamma},{\delta\over\delta
x^j}\right){\partial\over\partial x^i_\alpha}=P^{(l)(\alpha)\;(\gamma)}_
{(\delta)(i)j(k)}{\partial\over\partial x^l_\delta}},
$$
$$
\displaystyle{
\mbox{\em\bf R}\left({\partial\over\partial x^k_\gamma},{\partial\over\partial
x^j_\beta}\right){\delta\over\delta t^\alpha}=\bar S^{\delta(\beta)(\gamma)}
_{\alpha(j)(k)}{\delta\over\delta t^\delta}},\quad
\displaystyle{
\mbox{\em\bf R}\left({\partial\over\partial x^k_\gamma},{\partial\over\partial
x^j_\beta}\right){\delta\over\delta x^i}=S^{l(\beta)(\gamma)}_{i(j)(k)}
{\delta\over\delta x^l}},
$$
$$
\displaystyle{
\mbox{\em\bf R}\left({\partial\over\partial x^k_\gamma},{\partial\over\partial
x^j_\beta}\right){\partial\over\partial x^i_\alpha}=S^{(l)(\alpha)(\beta)(
\gamma)}_{(\delta)(i)(j)(k)}{\partial\over\partial x^l_\delta}},
$$
which we arrange in the following table
\begin{equation}
\begin{tabular}{|c|c|c|c|}
\hline
&$h_T$&$h_M$&$v$\\
\hline
$h_Th_T$&$\bar R^\delta_{\alpha\beta\gamma}$&$R^l_{i\beta\gamma}$&$
R^{(l)(\alpha)}_{(\delta)(i)\beta\gamma}$\\
\hline
$h_Mh_T$&$\bar R^\delta_{\alpha\beta k}$&$R^l_{i\beta k}$&$
R^{(l)(\alpha)}_{(\delta)(i)\beta k}$\\
\hline
$h_Mh_M$&$\bar R^\delta_{\alpha jk}$&$R^l_{ijk}$&$
R^{(l)(\alpha)}_{(\delta)(i)jk}$\\
\hline
$vh_T$&$\bar P^{\delta\;\;(\gamma)}_{\alpha\beta(k)}$&$P^{l\;\;(\gamma)}_
{i\beta(k)}$&$P^{(l)(\alpha)\;(\gamma)}_{(\delta)(i)\beta(k)}$\\
\hline
$vh_M$&$\bar P^{\delta\;\;(\gamma)}_{\alpha j(k)}$&$P^{l\;\;(\gamma)}_
{ij(k)}$&$P^{(l)(\alpha)\;(\gamma)}_{(\delta)(i)j(k)}$\\
\hline
$vv$&$\bar S^{\delta(\beta)(\gamma)}_{\alpha(j)(k)}$&$S^{l(\beta)(\gamma)}_
{i(j)(k)}$&$S^{(l)(\alpha)(\beta)(\gamma)}_{(\delta)(i)(j)(k)}$\\
\hline
\end{tabular}
\end{equation}
\end{th}

Moreover, using the properties of the d-tensor {\bf R} and the expressions of
local $T$-, $M$-horizontal and vertical covariant derivatives attached to the $\Gamma$-linear
connection $\nabla$, we derive the following local components of the  curvature
d-tensor,\medskip\linebreak
$(h_T)\hspace{6mm}
\left\{\begin{array}{rl}\medskip
1.&\displaystyle{\bar R^\delta_{\alpha\beta\gamma}={\delta\bar G^\delta_
{\alpha\beta}\over\delta t^\gamma}-{\delta\bar G^\delta_
{\alpha\gamma}\over\delta t^\beta}+\bar G^\mu_{\alpha\beta}\bar G^\delta_{\mu
\gamma}-\bar G^\mu_{\alpha\gamma}\bar G^\delta_{\mu\beta}+\bar C^{\delta(\mu)}
_{\alpha(m)}R^{(m)}_{(\mu)\beta\gamma}}\\\medskip
2.&\displaystyle{\bar R^\delta_{\alpha\beta k}={\delta\bar G^\delta_
{\alpha\beta}\over\delta x^k}-{\delta\bar L^\delta_
{\alpha k}\over\delta t^\beta}+\bar G^\mu_{\alpha\beta}\bar L^\delta_{\mu
k}-\bar L^\mu_{\alpha k}\bar G^\delta_{\mu\beta}+\bar C^{\delta(\mu)}
_{\alpha(m)}R^{(m)}_{(\mu)\beta k}}\\\medskip
3.&\displaystyle{\bar R^\delta_{\alpha jk}={\delta\bar L^\delta_
{\alpha j}\over\delta x^k}-{\delta\bar L^\delta_
{\alpha k}\over\delta x^j}+\bar L^\mu_{\alpha j}\bar L^\delta_{\mu
k}-\bar L^\mu_{\alpha k}\bar L^\delta_{\mu j}+\bar C^{\delta(\mu)}
_{\alpha(m)}R^{(m)}_{(\mu)jk}}\\\medskip
4.&\displaystyle{\bar P^{\delta\;\;(\gamma)}_{\alpha\beta(k)}={\partial
\bar G^\delta_{\alpha\beta}\over\partial x^k_\gamma}-\bar C^{\delta(\gamma)}_
{\alpha(k)/\beta}+\bar C^{\delta(\mu)}_{\alpha(m)}P^{(m)\;\;(\gamma)}_{(\mu)
\beta(k)}}\\\medskip
5.&\displaystyle{\bar P^{\delta\;\;(\gamma)}_{\alpha j(k)}={\partial
\bar L^\delta_{\alpha j}\over\partial x^k_\gamma}-\bar C^{\delta(\gamma)}_
{\alpha(k)\vert j}+\bar C^{\delta(\mu)}_{\alpha(m)}P^{(m)\;\;(\gamma)}_{(\mu)
j(k)}}\\\medskip
6.&\displaystyle{\bar S^{\delta(\beta)(\gamma)}_{\alpha(j)(k)}={\partial\bar
C^{\delta(\beta)}_{\alpha(j)}\over\partial x^k_\gamma}-{\partial\bar C^{\delta
(\gamma)}_{\alpha(k)}\over\partial x^j_\beta}+\bar C^{\mu(\beta)}_{\alpha(j)}
\bar C^{\delta(\gamma)}_{\mu(k)}-\bar C^{\mu(\gamma)}_{\alpha(k)}\bar C^{\delta
(\beta)}_{\mu(j)}},
\end{array}\right.\medskip
$
$(h_M)\hspace{5mm}
\left\{\begin{array}{rl}\medskip
7.&\displaystyle{R^l_{i\beta\gamma}={\delta G^l_{i\beta}\over\delta t^\gamma}
-{\delta G^l_{i\gamma}\over\delta t^\beta}+G^m_{i\beta}G^l_{m\gamma}-
G^m_{i\gamma}G^l_{m\beta}+C^{l(\mu)}_{i(m)}R^{(m)}_{(\mu)\beta\gamma}}\\\medskip
8.&\displaystyle{R^l_{i\beta k}={\delta G^l_{i\beta}\over\delta x^k}-{\delta
L^l_{ik}\over\delta t^\beta}+G^m_{i\beta}L^l_{mk}-L^m_{ik}G^l_{m\beta}+
C^{l(\mu)}_{i(m)}R^{(m)}_{(\mu)\beta k}}\\\medskip
9.&\displaystyle{R^l_{ijk}={\delta L^l_{ij}\over\delta x^k}-{\delta L^l_
{ik}\over\delta x^j}+L^m_{ij}L^l_{mk}-L^m_{ik}L^l_{mj}+C^{l(\mu)}_{i(m)}
R^{(m)}_{(\mu)jk}}\\\medskip
10.&\displaystyle{P^{l\;\;(\gamma)}_{i\beta(k)}={\partial G^l_{i\beta}\over
\partial x^k_\gamma}-C^{l(\gamma)}_{i(k)/\beta}+C^{l(\mu)}_{i(m)}P^{(m)\;\;
(\gamma)}_{(\mu)\beta(k)}}\\\medskip
11.&\displaystyle{P^{l\;\;(\gamma)}_{ij(k)}={\partial L^l_{ij}\over\partial
x^k_\gamma}-C^{l(\gamma)}_{i(k)\vert j}+C^{l(\mu)}_{i(m)}P^{(m)\;\;(\gamma)}_
{(\mu)j(k)}}\\\medskip
12.&\displaystyle{S^{l(\beta)(\gamma)}_{i(j)(k)}={\partial C^{l(\beta)}_
{i(j)}\over\partial x^k_\gamma}-{\partial C^{l(\gamma)}_{i(k)}\over\partial
x^j_\beta}+C^{m(\beta)}_{i(j)}C^{l(\gamma)}_{m(k)}-C^{m(\gamma)}_{i(k)}
C^{l(\beta)}_{m(j)}},\\
\end{array}\right.\medskip
$
$(v)\hspace{8mm}
\left\{\begin{array}{rl}\medskip
13.&\displaystyle{R^{(l)(\alpha)}_{(\delta)(i)\beta\gamma}={\delta G^{(l)
(\alpha)}_{(\delta)(i)\beta}\over\delta t^\gamma}
-{\delta G^{(l)(\alpha)}_{(\delta)(i)\gamma}\over\delta t^\beta}+G^{(m)(\alpha)}
_{(\mu)(i)\beta}G^{(l)(\mu)}_{(\delta)(m)\gamma}-}\\\medskip
&\mbox{\hspace{13mm}}-G^{(m)(\alpha)}_{(\mu)(i)\gamma}G^{(l)(\mu)}_{(\delta)(m)\beta}
+C^{(l)(\alpha)(\mu)}_{(\delta)(i)(m)}R^{(m)}_{(\mu)\beta\gamma}\\\medskip
14.&\displaystyle{R^{(l)(\alpha)}_{(\delta)(i)\beta k}={\delta G^{(l)(\alpha)}
_{(\delta)(i)\beta}\over\delta x^k}-{\delta L^{(l)(\alpha)}_{(\delta)(i)k}
\over\delta t^\beta}+G^{(m)(\alpha)}_{(\mu)(i)\beta}L^{(l)(\mu)}_{(\delta)(m)
k}-}\\\medskip
&\mbox{\hspace{13mm}}-L^{(m)(\alpha)}_{(\mu)(i)k}G^{(l)(\mu)}_{(\delta)(m)
\beta}+C^{(l)(\alpha)(\mu)}_{(\delta)(i)(m)}R^{(m)}_{(\mu)\beta k}\\\medskip
15.&\displaystyle{R^{(l)(\alpha)}_{(\delta)(i)jk}={\delta L^{(l)(\alpha)}
_{(\delta)(i)j}\over\delta x^k}-{\delta L^{(l)(\alpha)}_{(\delta)(i)k}
\over\delta x^j}+L^{(m)(\alpha)}_{(\mu)(i)j}L^{(l)(\mu)}_{(\delta)(m)k}-}\\
\medskip
&\mbox{\hspace{13mm}}-L^{(m)(\alpha)}_{(\mu)(i)k}L^{(l)(\mu)}_{(\delta)(m)j}+
C^{(l)(\alpha)(\mu)}_{(\delta)(i)(m)}R^{(m)}_{(\mu)jk}\\\medskip
16.&\displaystyle{P^{(l)(\alpha)\;\;(\gamma)}_{(\delta)(i)\beta(k)}={\partial
G^{(l)(\alpha)}_{(\delta)(i)\beta}\over\partial x^k_\gamma}-C^{(l)(\alpha)
(\gamma)}_{(\delta)(i)(k)/\beta}+C^{(l)(\alpha)(\mu)}_{(\delta)(i)(m)}
P^{(m)\;\;(\gamma)}_{(\mu)\beta(k)}}\\\medskip
17.&\displaystyle{P^{(l)(\alpha)\;\;(\gamma)}_{(\delta)(i)j(k)}={\partial
L^{(l)(\alpha)}_{(\delta)(i)j}\over\partial x^k_\gamma}-C^{(l)(\alpha)
(\gamma)}_{(\delta)(i)(k)\vert j}+C^{(l)(\alpha)(\mu)}_{(\delta)(i)(m)}
P^{(m)\;\;(\gamma)}_{(\mu)j(k)}}\\\medskip
18.&\displaystyle{S^{(l)(\alpha)(\beta)(\gamma)}_{(\delta)(i)(j)(k)}={\partial
C^{(l)(\alpha)(\beta)}_{(\delta)(i)(j)}\over\partial x^k_\gamma}-{\partial
C^{(l)(\alpha)(\gamma)}_{(\delta)(i)(k)}\over\partial x^j_\beta}+C^{(m)
(\alpha)(\beta)}_{(\mu)(i)(j)}C^{(l)(\mu)(\gamma)}_{(\delta)(m)(k)}-}\\\medskip
&\mbox{\hspace{18mm}}-C^{(m)(\alpha)(\gamma)}_{(\mu)(i)(k)}C^{(l)(\mu)(\beta)}
_{(\delta)(m)(j)}.\\
\end{array}\right.\medskip
$\medskip\linebreak
\addtocounter{rem}{1}
{\bf Remark \therem} In the case of the Berwald $\Gamma_0$-linear connection associated to
the metrics pair $(h_{\alpha\beta},\varphi_{ij})$, all curvature d-tensors
vanish, except
$$
R^\delta_{\alpha\beta\gamma}=H^\delta_{\alpha\beta\gamma},\quad
R^l_{ijk}=r^l_{ijk},
$$
where $H^\delta_{\alpha\beta\gamma}$ (resp. $r^l_{ijk}$) are the curvature
tensors of the metric $h_{\alpha\beta}$ (resp. $\varphi_{ij}$).

\section{Ricci and Bianchi identities}

\setcounter{equation}{0}
\hspace{5mm} Taking into account the local form of the $T$- ,$M$-horizontal and
vertical covariant derivatives defined in Section 1, by a direct calculation one proves
\begin{th}
If $\displaystyle{X=X^{\alpha}{\delta\over\delta t^\alpha}+
X^i{\delta\over\delta x^i}+X^{(i)}_{(\alpha)}{\partial\over\partial
x^i_\alpha}}$ is an arbitrary vector field on the 1-jet space $E$, then
the following Ricci identities hold good:\medskip\\
$(h_T)\hspace{6mm}
\left\{\begin{array}{l}\medskip
X^\alpha_{/\beta/\gamma}-X^\alpha_{/\gamma/\beta}=X^\mu\bar R^\alpha_{\mu\beta\gamma}
-X^\alpha_{/\mu}\bar T^\mu_{\beta\gamma}-X^\alpha\vert^{(\mu)}_{(m)}R^{(m)}_
{(\mu)\beta\gamma}\\\medskip
X^\alpha_{/\beta\vert k}-X^\alpha_{\vert k/\beta}=X^\mu\bar R^\alpha_
{\mu\beta k}-X^\alpha_{/\mu}\bar T^\mu_{\beta k}-X^\alpha_{\vert m}T^m_{\beta
k}-X^\alpha\vert^{(\mu)}_{(m)}R^{(m)}_{(\mu)\beta k}\\\medskip
X^\alpha_{\vert j\vert k}-X^\alpha_{\vert k\vert j}=X^\mu\bar R^\alpha_
{\mu jk}-X^\alpha_{\vert m}T^m_{jk}-X^\alpha\vert^{(\mu)}_{(m)}R^{(m)}_{(\mu)
jk}\\\medskip
X^\alpha_{/\beta}\vert^{(\gamma)}_{(k)}-X^\alpha\vert^{(\gamma)}_{(k)/\beta}=
X^\mu\bar P^{\alpha\;\;(\gamma)}_{\mu\beta(k)}-X^\alpha_{/\mu}\bar C^{\mu(\gamma)}
_{\beta(k)}-X^\alpha\vert^{(\mu)}_{(m)}P^{(m)\;\;(\gamma)}_{(\mu)\beta(k)}
\\\medskip
X^\alpha_{\vert j}\vert^{(\gamma)}_{(k)}-X^\alpha\vert^{(\gamma)}_{(k)\vert j}=
X^\mu\bar P^{\alpha\;\;(\gamma)}_{\mu j(k)}-X^\alpha_{\vert m}C^{m(\gamma)}
_{j(k)}-X^\alpha\vert^{(\mu)}_{(m)}P^{(m)\;\;(\gamma)}_{(\mu)j(k)}
\\\medskip
X^\alpha\vert_{(j)}^{(\beta)}\vert^{(\gamma)}_{(k)}-X^\alpha\vert^{(\gamma)}_
{(k)}\vert^{(\beta)}_{(j)}=X^\mu\bar S^{\alpha(\beta)(\gamma)}_{\mu(j)(k)}-
X^\alpha\vert^{(\mu)}_{(m)}S^{(m)(\beta)(\gamma)}_{(\mu)(j)(k)},
\end{array}\right.\medskip
$
$(h_M)\hspace{5mm}
\left\{\begin{array}{l}\medskip
X^i_{/\beta/\gamma}-X^i_{/\gamma/\beta}=X^m R^i_{m\beta\gamma}
-X^i_{/\mu}\bar T^\mu_{\beta\gamma}-X^i\vert^{(\mu)}_{(m)}R^{(m)}_
{(\mu)\beta\gamma}\\\medskip
X^i_{/\beta\vert k}-X^i_{\vert k/\beta}=X^m R^i_{m\beta k}-X^i_{/\mu}\bar
T^\mu_{\beta k}-X^i_{\vert m}T^m_{\beta k}-X^i\vert^{(\mu)}_{(m)}R^{(m)}_
{(\mu)\beta k}\\\medskip
X^i_{\vert j\vert k}-X^i_{\vert k\vert j}=X^m R^i_{mjk}-X^i_{\vert m}T^m_{jk}
-X^i\vert^{(\mu)}_{(m)}R^{(m)}_{(\mu)jk}\\\medskip
X^i_{/\beta}\vert^{(\gamma)}_{(k)}-X^i\vert^{(\gamma)}_{(k)/\beta}=
X^m P^{i\;\;(\gamma)}_{m\beta(k)}-X^i_{/\mu}\bar C^{\mu(\gamma)}
_{\beta(k)}-X^i\vert^{(\mu)}_{(m)}P^{(m)\;\;(\gamma)}_{(\mu)\beta(k)}
\\\medskip
X^i_{\vert j}\vert^{(\gamma)}_{(k)}-X^i\vert^{(\gamma)}_{(k)\vert j}=
X^m P^{i\;\;(\gamma)}_{mj(k)}-X^i_{\vert m}C^{m(\gamma)}
_{j(k)}-X^i\vert^{(\mu)}_{(m)}P^{(m)\;\;(\gamma)}_{(\mu)j(k)}
\\\medskip
X^i\vert_{(j)}^{(\beta)}\vert^{(\gamma)}_{(k)}-X^i\vert^{(\gamma)}_
{(k)}\vert^{(\beta)}_{(j)}=X^m S^{i(\beta)(\gamma)}_{m(j)(k)}-
X^i\vert^{(\mu)}_{(m)}S^{(m)(\beta)(\gamma)}_{(\mu)(j)(k)},
\end{array}\right.\medskip
$
$(v)\hspace{8mm}
\left\{\begin{array}{l}\medskip
X^{(i)}_{(\alpha)/\beta/\gamma}-X^{(i)}_{(\alpha)/\gamma/\beta}=X^{(m)}_{(\mu)}
R^{(i)(\mu)}_{(\alpha)(m)\beta\gamma}-X^{(i)}_{(\alpha)/\mu}\bar T^\mu_{\beta
\gamma}-X^{(i)}_{(\alpha)}\vert^{(\mu)}_{(m)}R^{(m)}_
{(\mu)\beta\gamma}\\\medskip
X^{(i)}_{(\alpha)/\beta\vert k}-X^{(i)}_{(\alpha)\vert k/\beta}=X^{(m)}_{(\mu)}
R^{(i)(\mu)}_{(\alpha)(m)\beta k}-X^{(i)}_{(\alpha)/\mu}\bar T^\mu_{\beta k}-
X^{(i)}_{(\alpha)\vert m}T^m_{\beta k}-\\\medskip
\mbox{\hspace{30mm}}-X^{(i)}_{(\alpha)}\vert^{(\mu)}_{(m)}R^{(m)}_{(\mu)\beta
k}\\\medskip
X^{(i)}_{(\alpha)\vert j\vert k}-X^{(i)}_{(\alpha)\vert k\vert j}=X^{(m)}_{(\mu)}
R^{(i)(\mu)}_{(\alpha)(m)jk}-X^{(i)}_{(\alpha)\vert m}T^m_{jk}-X^{(i)}_{(\alpha)}
\vert^{(\mu)}_{(m)}R^{(m)}_{(\mu)jk}\\\medskip
X^{(i)}_{(\alpha)/\beta}\vert^{(\gamma)}_{(k)}-X^{(i)}_{(\alpha)}\vert^{(
\gamma)}_{(k)/\beta}=X^{(m)}_{(\mu)} P^{(i)(\mu)\;\;(\gamma)}_{(\alpha)(m)
\beta(k)}-X^{(i)}_{(\alpha)/\mu}\bar C^{\mu(\gamma)}_{\beta(k)}-\\\medskip
\mbox{\hspace{35mm}}-X^{(i)}_{(\alpha)}\vert^{(\mu)}_{(m)}P^{(m)\;\;(\gamma)}
_{(\mu)\beta(k)}\\\medskip
X^{(i)}_{(\alpha)\vert j}\vert^{(\gamma)}_{(k)}-X^{(i)}_{(\alpha)}\vert^{(
\gamma)}_{(k)\vert j}=X^{(m)}_{(\mu)} P^{(i)(\mu)\;\;(\gamma)}_{(\alpha)(m)j
(k)}-X^{(i)}_{(\alpha)\vert m}C^{m(\gamma)}_{j(k)}-\\\medskip
\mbox{\hspace{33mm}}-X^{(i)}_{(\alpha)}\vert^{
(\mu)}_{(m)}P^{(m)\;\;(\gamma)}_{(\mu)j(k)}
\\\medskip
X^{(i)}_{(\alpha)}\vert_{(j)}^{(\beta)}\vert^{(\gamma)}_{(k)}-X^{(i)}_{(\alpha)}
\vert^{(\gamma)}_{(k)}\vert^{(\beta)}_{(j)}=X^{(m)}_{(\mu)} S^{(i)(\mu)(\beta)
(\gamma)}_{(\alpha)(m)(j)(k)}-X^{(i)}_{(\alpha)}\vert^{(\mu)}_{(m)}S^{(m)(
\beta)(\gamma)}_{(\mu)(j)(k)}.
\end{array}\right.\medskip
$
\end{th}
\addtocounter{rem}{1}
{\bf Remark \therem} For the arbitrary vector fields $X,Y,Z\in{\cal X}(E)$ and
the arbitrary 1-form $\omega\in{\cal X}^*(E)$ on $J^1(T,M)$, the relations
\begin{equation}
\left\{
\begin{array}{l}\medskip
\mbox{\bf R}(X,Y)\omega=-\omega\circ\mbox{\bf R}(X,Y),\\
\mbox{\bf R}(X,Y)(Z\otimes\omega)=\mbox{\bf R}(X,Y)Z\otimes\omega+Z\otimes
\mbox{\bf R}(X,Y)\omega,
\end{array}\right.
\end{equation}
are true. These relations allow us to generalize the Ricci identities to the d-tensors
set of the 1-jet fibre bundle $E$. The generalization is a natural one, but
the expressions of Ricci identities become extremely complicated. For that
reason, we exemplify this generalization writing just one Ricci identity.
For example, if
$D=(D^{\delta i(k)(\varepsilon)\ldots}_{\alpha j(\eta)(l)\ldots})$
is an arbitrary d-tensor field on $E$, then the following Ricci identity
$$
\begin{array}{l}\medskip
D^{\delta i(k)(\varepsilon)\ldots}_{\alpha j(\eta)(l)\ldots/\beta/\gamma}-
D^{\delta i(k)(\varepsilon)\ldots}_{\alpha j(\eta)(l)\ldots/\gamma/\beta}=
D^{\mu i(k)(\varepsilon)\ldots}_{\alpha j(\eta)(l)\ldots}\bar R^\alpha_{\mu
\beta\gamma}+D^{\delta m(k)(\varepsilon)\ldots}_{\alpha j(\eta)(l)\ldots}
R^i_{m\beta\gamma}+\\\medskip
+D^{\delta i(m)(\varepsilon)\ldots}_{\alpha j(\mu)(l)\ldots}R^{(k)(\mu)}_
{(\eta)(m)\beta\gamma}+\ldots\ldots\ldots-D^{\delta i(k)(\varepsilon)\ldots}_
{\mu j(\eta)(l)\ldots}\bar R^\mu_{\alpha\beta\gamma}-D^{\delta i(k)(
\varepsilon)\ldots}_{\alpha m(\eta)(l)\ldots}R^m_{i\beta\gamma}-\\\medskip
-D^{\delta i(k)(\mu)\ldots}_{\alpha j(\eta)(m)\ldots}R^{(m)(\varepsilon)}_
{(\mu)(l)\beta\gamma}-\ldots\ldots\ldots-D^{\delta i(k)(\varepsilon)\ldots}_
{\alpha j(\eta)(l)\ldots/\mu}\bar T^\mu_{\beta\gamma}-D^{\delta i(k)(
\varepsilon)\ldots}_{\alpha j(\eta)(l)\ldots}\vert^{(\mu)}_{(m)}R^{(m)}_{(\mu
)\beta\gamma}
\end{array}
$$
holds good.

Now, let us consider  the Liouville canonical vector field {\bf C}$\displaystyle
{=x^i_\alpha{\partial\over\partial x^i_\alpha}}$ and the {\it deflection
d-tensors associated to the $\Gamma$-linear connection $\nabla$}, defined by
the local components
$$
\bar D^{(i)}_{(\alpha)\beta}=x^i_{\alpha/\beta},\quad
D^{(i)}_{(\alpha)j}=x^i_{\alpha\vert j},\quad
d^{(i)(\beta)}_{(\alpha)(j)}=x^i_\alpha\vert^{(\beta)}_{(j)}.
$$
By a direct calculation, we find
\begin{equation}
\left\{\begin{array}{l}\medskip
\bar D^{(i)}_{(\alpha)\beta}=-M^{(i)}_{(\alpha)\beta}+G^{(i)(\mu)}_{(\alpha)(
m)\beta}x^m_\mu\\\medskip
D^{(i)}_{(\alpha)j}=-N^{(i)}_{(\alpha)j}+L^{(i)(\mu)}_{(\alpha)(m)j}
x^m_\mu\\\medskip
d^{(i)(\beta)}_{(\alpha)(j)}=\delta^i_j\delta^\beta_\alpha+C^{(i)(\mu)(
\beta)}_{(\alpha)(m)(j)}x^m_\mu.
\end{array}\right.
\end{equation}

Applying the $v$-set of the Ricci identities to the components of the Liouville
vector field, we obtain
\begin{th}
The deflection d-tensors, attached to the $\Gamma$-linear connection $\nabla$,
satisfy:
$$
\left\{\begin{array}{l}\medskip
\bar D^{(i)}_{(\alpha)\beta/\gamma}-\bar D^{(i)}_{(\alpha)\gamma/\beta}=
x^m_\mu R^{(i)(\mu)}_{(\alpha)(m)\beta\gamma}-\bar D^{(i)}_{(\alpha)\mu}
\bar T^\mu_{\beta\gamma}-d^{(i)(\mu)}_{(\alpha)(m)}R^{(m)}_{(\mu)\beta\gamma}
\\\medskip
\bar D^{(i)}_{(\alpha)\beta\vert k}-D^{(i)}_{(\alpha)k/\beta}=x^m_\mu
R^{(i)(\mu)}_{(\alpha)(m)\beta k}-\bar D^{(i)}_{(\alpha)\mu}\bar T^\mu_{\beta
k}-D^{(i)}_{(\alpha)m}T^m_{\beta k}-d^{(i)(\mu)}_{(\alpha)(m)}R^{(m)}_{(
\mu)\beta k}\\\medskip
D^{(i)}_{(\alpha)j\vert k}-D^{(i)}_{(\alpha)k\vert j}=x^m_\mu
R^{(i)(\mu)}_{(\alpha)(m)jk}-D^{(i)}_{(\alpha)m}T^m_{jk}-d^{(i)(\mu)}_{(
\alpha)(m)}R^{(m)}_{(\mu)jk}\\\medskip
\bar D^{(i)}_{(\alpha)\beta}\vert^{(\gamma)}_{(k)}-d^{(i)(\gamma)}_{(\alpha)
(k)/\beta}=x^m_\mu P^{(i)(\mu)\;\;(\gamma)}_{(\alpha)(m)\beta(k)}-\bar D^
{(i)}_{(\alpha)\mu}\bar C^{\mu(\gamma)}_{\beta(k)}-d^{(i)(\mu)}_{(\alpha)(m)}
P^{(m)\;\;(\gamma)}_{(\mu)\beta(k)}\\\medskip
D^{(i)}_{(\alpha)j}\vert^{(\gamma)}_{(k)}-d^{(i)(\gamma)}_{(\alpha)(k)\vert j}
=x^m_\mu P^{(i)(\mu)\;\;(\gamma)}_{(\alpha)(m)j(k)}-D^{(i)}_{(\alpha)m}
C^{m(\gamma)}_{j(k)}-d^{(i)(\mu)}_{(\alpha)(m)}P^{(m)\;\;(\gamma)}_{(\mu)j(k)}
\\\medskip
d^{(i)(\beta)}_{(\alpha)(j)}\vert_{(k)}^{(\gamma)}-d^{(i)(\gamma)}_{(\alpha)
(k)}\vert^{(\beta)}_{(j)}=x^m_\mu S^{(i)(\mu)(\beta)(\gamma)}_{(\alpha)(m)(j)
(k)}-d^{(i)(\mu)}_{(\alpha)(m)}S^{(m)(\beta)(\gamma)}_{(\mu)(j)(k)}.
\end{array}\right.
$$
\end{th}

Finally, note that the torsion {\bf T} and the curvature {\bf R} of the
$\Gamma$-linear connection $\nabla$ are not independent. They verify the
general Bianchi identities \cite{5}
$$\left\{
\begin{array}{l}\medskip
\sum_{\{X,Y,Z\}}\{(\nabla_X\mbox{\bf T})(Y,Z)-\mbox{\bf R}(X,Y)Z+\mbox{\bf
T}(\mbox{\bf T}(X,Y),Z)\}=0,\quad\forall\;X,Y,Z\in{\cal X}(E)\\
\sum_{\{X,Y,Z\}}\{(\nabla_X\mbox{\bf R})(U,Y,Z)+\mbox{\bf R}(\mbox{\bf T}(X,Y
),Z)U\}=0,\quad\forall\;X,Y,Z,U\in{\cal X}(E),
\end{array}\right.\nopagebreak
$$
where $\{X,Y,Z\}$ means cyclic sum.

In the adapted basis $\displaystyle{(X_A)=\left({\delta\over\delta t^\alpha},
{\delta\over\delta x^i},{\partial\over\partial x^i_\alpha}\right)}$, we have
sixty-four effective  Bianchi identities, obtained by the relations
\begin{equation}
\left\{\begin{array}{l}\medskip
\sum_{\{A,B,C\}}\{R^F_{ABC}-T^F_{AB:C}-T^G_{AB}T^F_{CG}\}=0\\
\sum_{\{A,B,C\}}\{R^F_{DAB:C}+T^G_{AB}R^F_{DAG}\}=0,
\end{array}\right.
\end{equation}
where {\bf R}$(X_A,X_B)X_C=R^D_{CBA}X_D$, {\bf T}$(X_A,X_B)=T^D_{BA}X_D$ and
"$_{:A}$" represents one of the covariant derivatives $"_{/\alpha}"$, $"_{\vert i}"$
or $"\vert^{(\alpha)}_{(i)}"$. The large number and the complicated form of the
above Bianchi identities associated to a $\Gamma$-linear connection determine
us to study them, in a subsequent paper \cite{9}, in the more particular case
of the {\it $h$-normal $\Gamma$-linear connections}. In that case, the number
of Bianchi identities reduces to thirty.

\section{Jet prolongation of vector fields}

\setcounter{equation}{0}
\hspace{5mm} A general vector field $X^*$ on $J^1(T,M)$ can be written under
the form
$$
X^*=X^\alpha{\partial\over\partial t^\alpha}+X^i{\partial\over\partial x^i}
+X^{(i)}_{(\alpha)}{\partial\over\partial x^i_\alpha},
$$
where the components $X^\alpha,\;X^i\;X^{(i)}_{(\alpha)}$ are functions of $(t^\alpha,x^i,x^i_\alpha)$.

The prolongation of a vector field $X$ on $T\times M$ to a vector field on the
1-jet bundle $J^1(T,M)$ was solved by Olver \cite{12} in the following sense.
\medskip\\
\addtocounter{defin}{1}
{\bf Definition \thedefin} Let $X$ be a vector field on $T\times M$ with corresponding
(local) one-parameter group $\exp(\varepsilon X)$. The {\it 1-th prolongation}
of $X$, denoted by $pr^{(1)}X$, will be a vector field on the 1-jet space
$J^1(T,M)$, and is defined to be the infinitesimal generator of the corresponding
prolonged one-parameter group $pr^{(1)}[\exp(\varepsilon X)]$, i. e. ,
\begin{equation}
\displaystyle{[pr^{(1)}X](t^\alpha,x^i,x^i_\alpha)=\left.{d\over d\varepsilon}
\right\vert_{\varepsilon=0}pr^{(1)}[\exp(\varepsilon X)](t^\alpha,x^i,x^i_
\alpha).}
\end{equation}

In order to write the components of the prolongation, Olver used the {\it
$\alpha$-th total derivative} $D_\alpha$ of an arbitrary
function $f(t^\alpha,x^i)$ on $T\times M$, which is defined by the relation
\begin{equation}
\displaystyle{D_\alpha f={\partial f\over\partial t^\alpha}+{\partial f\over
\partial x^i}x^i_\alpha.}
\end{equation}

Thus, starting with $\displaystyle{X=X^\alpha(t,x){\partial\over\partial t^\alpha}+
X^i(t,x){\partial\over\partial x^i}}$ like a vector field on $T\times M$,
Olver introduced the 1-th prolongation of $X$ as the vector field
\begin{equation}
pr^{(1)}X=X+X^{(i)}_{(\alpha)}(t^\beta,x^j,x^j_\beta){\partial
\over\partial x^i_\alpha},
\end{equation}
where
$$
X^{(i)}_{(\alpha)}=D_\alpha X^i-(D_\alpha X^\beta)x^i_\beta={\partial X^i\over
\partial t^\alpha}+{\partial X^i\over\partial x^j}x^j_\alpha-\left({\partial
X^\beta\over\partial t^\alpha}+{\partial X^\beta\over\partial x^j}x^j_\alpha
\right)x^i_\beta.
$$

Let us use a geometrical approach for obtaining jet prolongations of vector fields.
If we assume that is given a nonlinear connection $\Gamma=(M^{(i)}_{(\alpha)\beta},
N^{(i)}_{(\alpha)j})$ on $J^1(T,M)$, then the $\alpha$-th total
derivative used by Olver can be written as
\begin{equation}
D_\alpha f={\delta f\over\delta t^\alpha}+{\delta f\over\delta x^i}x^i_\alpha
=f_{/\alpha}+f_{\vert i}x^i_\alpha,
\end{equation}
and, consequently, $D_\alpha f$ represent the local components
of a distinguished 1-form on $J^1(T,M)$, which is expressed by $Df=(D_\alpha
f)dt^\alpha$.

Now, let there be given a vector field $X$ on $T\times M$. If we suppose that
$J^1(T,M)$ is endowed at the same time with a $\Gamma$-linear connection
\ref{lglc}, we can define the geometrical 1-th jet prolongation of $X$,
\begin{equation}
pr^{(1)}X=X^\alpha{\delta\over\delta t^\alpha}+X^i{\delta\over\delta x^i}+
Y^{(i)}_{(\alpha)}(t^\beta,x^j,x^j_\beta){\partial\over\partial x^i_\alpha},
\end{equation}
setting
$$
Y^{(i)}_{(\alpha)}=X^i_{/\alpha}+X^i_{\vert j}x^j_\alpha-X^\beta_{/\alpha}
x^i_\beta-X^\beta_{\vert j}x^j_\alpha x^i_\beta+
$$
$$
\hspace*{7mm}+X^\mu\left(M^{(i)}_{(\alpha)\mu}+\bar G^\beta_{\mu\alpha}x^i_\beta+\bar
L^\beta_{\mu j}x^j_\alpha x^i_\beta\right)-
$$
$$
-X^m\left(N^{(i)}_{(\alpha)m}+G^i_{m\alpha}+L^i_{mj}x^j_\alpha\right).
$$
\addtocounter{rem}{1}
{\bf Remarks \therem} i) Our prolongation coincides with that of Olver. Moreover,
we have the relation
\begin{equation}
Y^{(i)}_{(\alpha)}=X^{(i)}_{(\alpha)}+M^{(i)}_{(\alpha)\mu}X^\mu+N^{(i)}_
{(\alpha)m}X^m.
\end{equation}

ii) In the particular case of the Berwald $\Gamma_0$-linear connection associated to the
metrics $h_{\alpha\beta}$ and $\varphi_{ij}$, the expression of $Y^{(i)}_{(
\alpha)}$ reduces to
\begin{equation}
Y^{(i)}_{(\alpha)}=X^i_{/\!/\alpha}+X^i_{\Vert j}x^j_\alpha-X^\beta_{/\!/\alpha}
x^i_\beta-X^\beta_{\Vert j}x^j_\alpha x^i_\beta-2\gamma^i_{jm}x^j_\alpha X^m,
\end{equation}
where $\gamma^i_{jk}$ represent the Christoffel symbols of the metric
$\varphi_{ij}$.\medskip\\
{\bf \underline{Open problem}.}

Study the prolongations of vectors, 1-forms, tensors, $G$-structures from
$T\times M$ to $J^1(T,M)$.\medskip\\
{\bf Acknowledgments.} It is a pleasure for us to thank Prof. Dr. D. Opri\c s
for many helpful discussions on this research. We also thank to the reviewers
of Kodai Mathematical Journal for their valuable comments upon a previous version
of this paper.
\pagebreak

\begin{center}
University POLITEHNICA of Bucharest\\
Department of Mathematics I\\
Splaiul Independentei 313\\
77206 Bucharest, Romania\\
e-mail: mircea@mathem.pub.ro\\
e-mail: udriste@mathem.pub.ro
\end{center}

\end{document}